\documentstyle[12pt]{article}
\newlength{\abstractwidth}
\renewcommand{\thefootnote}{\fnsymbol{footnote}}
\renewcommand{\thanks}[1]{\footnote{#1}} % Use this for footnotes
\newcommand{\starttext}{ \setcounter{footnote}{0}
\renewcommand{\thefootnote}{\arabic{footnote}}}

\newcommand{\be}{\begin{equation}}
\newcommand{\bea}{\begin{eqnarray}}
\newcommand{\eea}{\end{eqnarray}} \newcommand{\ee}{\end{equation}}
 \newcommand{\<}{\langle}
\renewcommand{\>}{\rangle} \def\ba{\begin{eqnarray}}
\def\ea{\end{eqnarray}}

\def\o{\omega}

\def\det{{\rm det}}

\def\log{\,{\rm log}\,}

\def\o{\omega}

\def\o{\omega}

\def\na{\nabla}

\def\p{\partial}

\def\ddb{{\partial\bar\partial}}

\def\na{{\nabla}}

\def\[{{\bf [}}
\def\]{{\bf ]}}

\begin{document}
\starttext \baselineskip=18pt \setcounter{footnote}{0}
\newtheorem{theorem}{Theorem}
\newtheorem{lemma}{Lemma}
\newtheorem{corollary}{Corollary}
\newtheorem{definition}{Definition}
\newtheorem{conjecture}{Conjecture}
\newtheorem{proposition}{Proposition}

\begin{center}
{\Large \bf GEOMETRIC FLOWS AND STROMINGER SYSTEMS
\footnote{Work supported in part by the National Science Foundation under Grant DMS-12-66033 and DMS-1308136.
Keywords: balanced metrics, anomaly equation, Nash-Moser implicit function theorem. 
%AMS classification numbers: 32Q26 (32Q15, 32Q20, 32U05, 32W20), 35Kxx.
}}

\bigskip

{\large Duong H. Phong, Sebastien Picard and Xiangwen Zhang} \\

\medskip

\begin{abstract}

\medskip
\small{
A geometric flow on $(2,2)$-forms is introduced which preserves the balanced condition of metrics, and whose stationary points satisfy the anomaly equation in Strominger systems. The existence of solutions for a short time is established, using Hamilton's version of the Nash-Moser implicit function theorem.
}

\end{abstract}

\end{center}

\baselineskip=15pt
\setcounter{equation}{0}
\setcounter{footnote}{0}

\section{Introduction}
\setcounter{equation}{0}

The Strominger system is a system of equations for a metric on a $3$-dimensional complex manifold $X$ equipped with a nowhere vanishing holomorphic $3$-form, and a Hermitian metric on a holomorphic vector bundle over $X$. It is of considerable interest both in physics, where it is the equation for supersymmetric compactifications of the heterotic string to a four-dimensional space-time, and in mathematics, where it is a non-K\"ahler generalization of a Calabi-Yau metric \cite{Y} coupled to a Hermitian-Einstein connection \cite{D, UY}. The first mathematically rigorous solutions to the Strominger system were found by Li-Yau \cite{LY} and Fu-Yau \cite{FY,FY2}, and more solutions were found in \cite{AG1,AG2,Fe1,Fe2,FeY,FIUV1,FIUV2,FTY,Gr,OUV,UV}. Other solutions on physical grounds were studied by physicists in e.g. \cite{BBFT,CI,DRS}. Generalizations of the Fu-Yau equation to higher dimensions are considered in
\cite{PPZ,PPZ0}. The general solution appears out of reach at the present time.

\smallskip
The main goal of this paper is to propose a geometric flow of $(2,2)$-forms, whose stationary points provide solutions of the Strominger system. We call it the Anomaly flow, as its curvature terms are the local characteristic classes arising in gravitational and Yang-Mills anomalies in string theory.
In this paper, we present some evidence that the Anomaly flow may provide a viable approach to Strominger systems. In particular, we show that it preserves the balanced property of metrics, that it also admits a more conventional description as a flow of Hermitian metrics by their curvatures, and that short-time solutions always exist for small values of the string tension parameter. 
Estimates for solutions and criteria for long-time existence and convergence are relegated to later work.

\section{The Anomaly flow}

Let $X$ be a compact $3$-dimensional complex manifold, which admits a nowhere vanishing holomorphic $(3,0)$-form $\Omega$. Let $E\to X$ be a holomorphic vector bundle over $X$. Let $\omega_0$ be a Hermitian metric on $X$, and $H_0$ a Hermitian metric on $E$. We define the Anomaly flow for the pair $(X,E)$
to be the flow $\o(t), H(t)$ of metrics on $X$ and on $E$ given by,
\bea
\label{anomalyXE}
\p_t (\|\Omega\|_\o\o^2 )
&=&
i\ddb\o-{\alpha' \over 4}\left({\rm Tr}(R\wedge R)-{\rm Tr}(F\wedge F)\right)
\nonumber\\
H^{-1} \, \p_tH
&=&
-\Lambda F
\eea
with initial condition $\o(0)=\o_0$, $H(0)=H_0$. Here $\alpha'$ is a fixed positive parameter, called the string tension in the physics literature. The expressions $R=(R^p{}_q)$ and $F= (F^{\alpha}{}_{\beta})$ denote respectively the curvature of the Chern connections defined by $\o(t)$ and $H(t)$ on $X$ and on $E$. They are $(1,1)$-forms, valued respectively in the bundles $End(T^{1,0})$ and $End(E)$ of endomorphisms of $T^{1,0}$ and $E$. Our conventions are
\bea
[\na_j,\na_{\bar k}]V^p=R_{\bar kj}{}^p{}_q V^q,
\qquad
[\na_j,\na_{\bar k}]\varphi^\alpha=F_{\bar kj}{}^\alpha{}_\beta \varphi^\beta,
\eea
where $V^p$ and $\varphi^\alpha$ are respectively sections of $T^{1,0}$ and $E$ in a local trivialization. The Hermitian form $\omega$ is defined by $\omega= i \, g_{\bar k j} dz^j\wedge d\bar z^k$, and a $(p, q)$ form $\eta$ has component $\eta_{\bar k_1 \cdots \bar k_q j_1 \cdots j_p}$ given by
\bea\label{convention}
 \eta= {1\over p! q!} \sum \eta_{\bar k_1 \cdots \bar k_q j_1 \cdots j_p}\, dz^{j_p} \wedge \cdots \wedge dz^{j_1} \wedge d \bar z^{k_q} \wedge \cdots d\bar z ^{k_1}.
\eea
With this convention, the $(1,1)$-forms $R=(R^p{}_q)$ and $F= (F^{\alpha}{}_{\beta})$ are then given by 
\begin{eqnarray*}
R^p{}_q=  R_{\bar kj}{}^p{}_q\,dz^j\wedge d\bar z^k, \ \  \ F^{\alpha}{}_{\beta}= F_{\bar kj}{}^\alpha{}_\beta\,dz^j\wedge d\bar z^k.
\end{eqnarray*}
We also define the pointwise inner product $\langle\, , \, \rangle_{\o}$ as
\bea\nonumber
\langle \phi, \psi \rangle_{\o} = {1\over p!q!} g^{\mu_1 \bar\beta_1} \cdots g^{\mu_q \bar \beta_q} g^{\alpha_1 \bar\lambda_1} \cdots g^{\alpha_p \bar \lambda_p}  \, \phi_{\bar\beta_1 \cdots \bar\beta_q \alpha_1\cdots \alpha_p}\overline{\psi_{\bar \mu_1 \cdots \bar\mu_q \lambda_1 \cdots \lambda_p}},
\eea
for any $(p, q)$ forms $\phi$ and $\psi$. The Hodge operator $\Lambda$ is defined as usual by
$(\Lambda F)^\alpha{}_\beta=g^{j\bar k}F_{\bar kj}{}^\alpha{}_\beta$. Here $(g^{j\bar k})$ denotes the inverse $(g_{\bar j k})^{-1}$ of $(g_{\bar j k})$.

\smallskip

We shall also be interested in a version of the Anomaly flow for just a metric $\o(t)$ on $X$. Thus let $\Phi_0$ be a given closed smooth $(2,2)$-form on $X$. The Anomaly flow for $X$ with given $\Phi_0$ is defined by
\bea
\label{anomalyX}
\p_t(\|\Omega\|_\o\o^2)
=
i\ddb\o-{\alpha' \over 4}({\rm Tr}\,R\wedge R-\Phi_0)
\eea
with initial condition $\o(0)=\o_0$. 

\smallskip
The following simple theorem provides the motivation for the Anomaly flows:

\begin{theorem}
\label{Th1}
Let $X,E,\Omega$ be as above, and consider the equations
(\ref{anomalyXE}) or (\ref{anomalyX}) on either $(X,E)$ or $X$, 
with initial metrics $\o_0$ and $H_0$.

\smallskip
{\rm (a)} The equation (\ref{anomalyX}) is well-defined as a flow, i.e., it defines a vector field on the space of metrics, or equivalently, a vector field on the space of positive Hermitian $(2,2)$-forms. The flows are local, in the sense that the vector fields are given by local expressions in the underlying Hermitian metrics (or the underlying $(2,2)$-form). 

Similarly for the equation (\ref{anomalyXE}), which defines now a vector field on the direct sum of the space of Hermitian metrics on $X$ (or positive $(2,2)$-forms) with the space of Hermitian metrics on $E$.

\smallskip
{\rm (b)} Assume that the Anomaly flows admit a smooth solution on some time interval $[0,T)$. If the initial metric $\o_0$ satisfies the balanced condition
\footnote{The usual balanced condition for a Hermitian metric $\o$ in dimension $n$ is $d\o^{n-1}=0$. For simplicity, we use the same terminology for the slight modification used in (\ref{balanced}).}
\bea
\label{balanced}
d(\|\Omega\|_{\o_0}\o_0^2)=0
\eea
then the metric $\o(t)$ will satisfy the same balanced condition, namely $d(\|\Omega\|_{\o_t}\o_t^2)=0$, for any $t$ in its time interval $[0,T)$ of existence. 

\smallskip
{\rm (c)} Assume that the Anomaly flow for $(X,E)$ exists for all time $[0,\infty)$, and that the initial metric satisfies the balanced-like condition (\ref{balanced}). If the flow $(\o(t),H(t))$ converges to a pair $(\o_\infty,H_\infty)$ of Hermitian metrics on $X$ and $E$, then this pair satisfies the Strominger system of equations
\footnote{It has recently been brought to our attention that the same system of equations for supersymmetric compactifications was independently proposed by
C. Hull \cite{Hull1,Hull2}.}
\bea
\label{S}
&&
F_\infty\wedge \omega_\infty^2=0
\nonumber\\
&&
F_\infty^{2,0}=F_\infty^{0,2}=0
\nonumber\\
&&
i\ddb\o_\infty={\alpha' \over 4}
\left({\rm Tr} R_\infty\wedge R_\infty-{\rm Tr} F_\infty\wedge F_\infty   \right)
\nonumber\\
&&
d\left(\|\Omega\|_{\o_\infty}\o_\infty^2\right)=0.
\eea
Similarly,
if the Anomaly flow on $X$ with given $\Phi_0$ converges, then the limiting metric $\o_\infty$ satisfies
\bea
&&
\label{SX}
i\ddb\o_\infty= {\alpha' \over 4}
({\rm Tr}(R_\infty\wedge R_\infty)-\Phi_0)
\nonumber\\
&&
d\left(\|\Omega\|_{\o_{\infty}}\o_{\infty}^2\right)=0.
\eea
\end{theorem}

\medskip

We note that the condition $d(\|\Omega\|_{\o_\infty}\o_\infty^2)=0$ has been shown by Li and Yau \cite{LY} to be equivalent to the condition $d^\dagger\omega_\infty=i(\bar\p-\p)\log \|\Omega\|_{\o_\infty}$, so that the system (\ref{S}) is indeed equivalent to the system originally written down by Strominger \cite{S}.

We also note that the equations for $H$ in the Strominger system are just the Hermitian-Einstein equation for a Chern unitary connection, and the flow of $H$ in (\ref{anomalyXE}) is of course the well-known Donaldson heat flow, which is gauge equivalent to the Yang-Mills flow.

\bigskip
Next, we consider the issue of short-time solutions for the flows, and when they would be parabolic. For a fixed metric $\o$, we define a modified operator $\tilde\star$ of the Hodge $\star_{\o}$ operator as the operator from the space of $(2,2)$-forms $\delta\Psi$ to the space of $(1,1)$-forms,
given by
\bea
\label{tildestar}
\tilde\star \,\delta\Psi
=
{1\over 2\|\Omega\|_\o}(\<\star \delta\Psi,\o\>\o-\star \delta\Psi).
\eea
Next, we view the curvature tensor $R_{\bar kj}{}^p{}_q$ of $\o$ as the operator $Rm$ from the space of $(1,1)$-forms $\delta\o$ into itself
given by
\bea
Rm(\delta\o)_{\bar kj}= \, R_{\bar kj}{}^{p\bar q}(\delta \o)_{\bar q p}.
\eea
We can then define the following linear differential operator $\tilde \Delta$ of order $2$ on the space of Hermitian tensors $\delta\Psi$ of type $(2,2)$,
\bea
\label{Delta}
\tilde\Delta\,(\delta\Psi)
=
i\ddb (\tilde\star \delta\Psi -{\alpha' \over 2} \,Rm(\tilde\star\delta\Psi)).
\eea
We note that the range of $\tilde\Delta$ is contained in the space of closed $(2,2)$-forms. We shall say that the operator $\tilde\Delta$ is {\it elliptic on the space of closed $(2,2)$-forms} if its symbol, restricted to the null space of the symbol of the exterior derivative $d$, admits only eigenvalues with strictly positive real parts.
We have then the following theorem:

\begin{theorem}
\label{Th2}
Let $(X,\Omega,E)$ be as before, and consider the flows (\ref{anomalyXE}) and (\ref{anomalyX}), with initial Hermitian metrics $\o_0$ and $H_0$ on $X$ and $E$ respectively. If the operator $\tilde\Delta$ with respect to the initial metric $\o_0$ is elliptic on the space of closed $(2,2)$-forms, in the above sense, then both flows (\ref{anomalyXE}) and (\ref{anomalyX}) admit a smooth solution on some non-trivial finite time interval.

\end{theorem}

\medskip

\section{Proof of Theorem \ref{Th1}}
\setcounter{equation}{0}

\par
Let $\Psi=\|\Omega\|_\o\o^2$. The essence of Theorem \ref{Th1} is that $\o$ can be recaptured from $\Psi$ (and of course vice versa), by purely local expressions. For this we need to discuss the issue of $(n-1)$-th root of a positive $(n-1, n-1)$-form in some detail.

\subsection{The $(n-1)$-th root of an $(n-1,n-1)$-form}

In general, in dimension $n$, let $\Phi$ be a $(n-1,n-1)$-form which is positive definite, in the sense that 
\begin{eqnarray*}
\Phi\wedge i \, \eta\wedge \bar\eta
\end{eqnarray*}
 is a positive $(n,n)$-form for any non-zero $(1,0)$-form $\eta$ and which equals $0$ if and only if $\eta=0$. Michelsohn \cite{M} has shown that there exists a unique positive $(1,1)$-form $\o$ with 
 \bea
 \o^{n-1}=\Phi.
 \eea
  We need a viable formula for $\o$, which can be obtained as follows. Let $\Phi$ be expressed as in
 \cite{M, TW1} by
\bea \label{n-1_convention}
 \Phi &=&   i^{n-1} (n-1)!\sum_{k, j} \, (sgn(k, j))\,\Phi^{k\bar j}dz^1\wedge d\bar z^1\wedge \cdots \wedge \widehat{dz^k}\wedge d\bar z^k\wedge \cdots \\\nonumber
 &&\qquad\qquad\qquad\wedge\, dz^j \wedge \widehat{d\bar z^j} \wedge \cdots\wedge dz^n \wedge d\bar z^n
\eea
where $sgn(k, j) = -1$ if $k>j$ and $sgn(k, j) = 1$ otherwise. One advantage for this representation is that $\Phi^{k\bar j}$ is a Hermitian matrix. Then the $(n-1)$-th root $\o= i \, g_{\bar j k} dz^k\wedge d\bar z^j$ of $\Phi$ is given by
\bea
\label{oPsi}
g_{\bar jk}=({\rm det}\,g)\,(\Phi^{-1})_{\bar jk},
\eea
where $(\Phi^{-1})_{\bar jk}$ is the inverse matrix of $\Phi^{k\bar j}$, i.e., 
$\Phi^{k\bar j}(\Phi^{-1})_{\bar j\ell}=\delta^k{}_\ell$. To see this, we note that the entry $(\o^{n-1})^{j\bar k}$ in the product $\o^{n-1}$ is obtained by taking the product of the entries, with corresponding permutation signs and $(n-1)!$ factor, of the matrix obtained from $g_{\bar pq}$ by removing the $j$-th row and the $k$-th column. In other words, it is the $j\bar k$ cofactor of the matrix $(g_{\bar pq})$. The equation (\ref{oPsi}) is just a reformulation of this statement.

\medskip
The notion of $(n-1)$-th root is independent of any metric. Nevertheless, it can be useful to express it in terms of the Hodge star operator.
If $\tilde\o= i\,  \tilde g_{\bar k j} dz^j \wedge d\bar z^k$ is any metric, recall that the Hodge star operator $\star_{\tilde\o}$ with respect to $\tilde\o$ is defined by the equation
\bea
\phi\wedge \overline\Phi={\<\phi, (\star_{\tilde\o}\Phi)\>_{\tilde\o}\over n!}\tilde\o^n
\eea
for any $(1,1)$-form $\phi= \phi_{\bar k j} dz^j \wedge d\bar z^k$. The left-hand side can be easily recognized to be
\bea
\label{star0}
i ^{-1}\, (n-1)! \, \phi_{\bar kj}\Phi^{j\bar k}
\prod_{\ell=1}^nidz^\ell\wedge d\bar z^{\ell}
=
- i\, {\phi_{\bar kj}\Phi^{j\bar k}\over {\rm det}\,\tilde g}
\,
{\tilde\o^n \over n}.
\eea
Here we use the fact that $(\Phi^{j\bar k})$ is Hermitian. We can also write the right-hand side as
\bea
 \tilde g^{p \bar k}\, \tilde g^{j \bar q}\, \phi_{\bar k j} \overline{(\star_{\tilde\o}\Phi)_{\bar p q}}\, {\tilde\o^n \over n!}.
\eea
This implies that 
\bea
\label{star}
(\star_{\tilde\o}\Phi)_{\bar pq}
=
i \, {(n-1)!\over {\rm det}\,\tilde g}\Phi^{k\bar j}\tilde g_{\bar p k}\tilde g_{\bar jq}.
\eea
As a check, we note that, since the expression in (\ref{star0}) is an $(n,n)$-form, the factor $\phi_{\bar kj}\Phi^{j\bar k}/{\rm det}\,\tilde g$ must be a scalar, and hence $\Phi^{k\bar j}$ should be interpreted as a section of $(\Lambda^{1,1})^*\otimes K_X\otimes \overline{K_X}$. This is consistent with the fact that the expression given in (\ref{star}) is a $(1,1)$-form.
Taking $\tilde\o$ in (\ref{star}) to be $\o$ itself, we obtain
\begin{eqnarray*}
\star_\o\Phi=(n-1)! \, \o,
\end{eqnarray*}
 a formula that can also be easily seen using an orthonormal basis for $\o$. Henceforth, we shall suppress the subindex $\o$ in the star operator, when the metric $\o$ is implicit.

\subsection{The relation $\Psi=\|\Omega\|_\o\o^2$}

We return to the setting of a $3$-fold $X$, equipped with a fixed nowhere vanishing $(3,0)$-form $\Omega$. We will use the notation $| \Omega |^2 = \Omega \overline{\Omega}$ and $\|\Omega\|_\o^2 = \Omega \overline{\Omega} (\det \, g)^{-1}$. If $\Psi$ is any positive $(2,2)$-form, we claim that there is a unique positive $(1,1)$-form $\o$ so that $\Psi=\|\Omega\|_\o\o^2$. 
\medskip
\par Indeed, this equation determines the norm $\|\Omega\|_\o$, since taking determinants gives
\bea
{\rm det}\,\Psi
=
\left(| \Omega |^2({\rm det}\,g)^{-1}\right)^{3/2}\,({\rm det}\,g)^2
\eea
and hence
\bea
({\rm det}\,g)^{1\over 2}={{\rm det}\,\Psi\over |\Omega |^3}.
\eea
This determines ${\rm det}\,g$ in terms of $\Psi$, and hence $\|\Omega\|_\o$ in terms of $\Psi$.  We can then obtain $\o$ as the square root of the positive $(2,2)$-form $\|\Omega\|_\o^{-1}\Psi$. The relations (\ref{oPsi}) and (\ref{star}) become
\bea
\label{root1}
g_{\bar jk}=
{{\rm det}\,\Psi\over |\Omega |^2}
\Psi_{\bar jk}^{-1},
\qquad
\star_\o
\Psi
=
2 \|\Omega\|_\o\o.
\eea

\subsection{Proof of Theorem \ref{Th1}, Part (a)}

It is now straightforward to relate the variations of $\o$ to the variations of $\Psi$. 
Differentiating the first equation on the left side of (\ref{root1}) gives
\bea
\label{variations}
\delta g_{\bar jk}
=
{1\over \|\Omega\|_\o\det g}
\left(g_{\bar qp}\, \delta\Psi^{p\bar q}\, g_{\bar jk}
-
g_{\bar jp}\,\delta\Psi^{p\bar q}\, g_{\bar qk}\right).
\eea
In intrinsic notation, using (\ref{star}), this can be rewritten as
\bea
\label{variations1}
\delta\o
=
{1\over 2 \|\Omega\|_\o}\left(\<\star\delta\Psi,\o\>\ \o-\star \delta\Psi\right)
=
\tilde\star\,\delta\Psi,
\eea
in view of the definition of the operator $\tilde\star$.

\smallskip
In particular, along the Anomaly flow, we can replace $\delta\o$ by $\p_t\o$ and $\delta\Psi$ by $\p_t\Psi=\p_t(\|\Omega\|_\o\o^2)$. We obtain in this way an equation giving $\p_t\o$ in terms of $\o$ and its curvature, which is the more conventional description of a geometric flow of Hermitian metrics. Equivalently, the flows can be written as flows of $(2,2)$-forms $\Psi$, given by a local vector field on the space of positive Hermitian $(2,2)$-forms.

\subsection{Proof of Theorem \ref{Th1}, Parts (b) and (c)}

The only non-trivial part of Theorem \ref{Th1} is the conceptual part due to the issue of taking square roots. Once this issue has been clarified, the proof is straightforward.

Part (b) follows immediately from the fact that the right hand sides of (\ref{anomalyXE}) and (\ref{anomalyX}) are always closed forms. This follows itself from the fact that $d\ddb\o=d^2\bar\partial\o=0$, and that both ${\rm Tr}(R\wedge R)$ and ${\rm Tr}(F\wedge F)$ are well-known closed representatives of the Chern classes $c_2(T^{1,0})$ and $c_2(E)$ of the bundles $T^{1,0}(X)$ and $E$. Thus
\bea
\p_t(d(\|\Omega\|_\o\o^2))
=
d\p_t(\|\Omega\|_\o\o^2)=0,
\eea
and $d(\|\Omega\|_\o\o^2)=0$ for all time $t$, if $d(\|\Omega\|_\o\o^2)=0$ at $t=0$.

\smallskip
Part (c) follows immediately from Part (b), which guarantees that the balanced-like equation in the Strominger systems is satisfied for all time. The equation $F_\infty^{2,0}=F_\infty^{0,2}=0$ is automatic for all Chern connections. The other equations are obvious consequences of the fact that the limiting metric $\o_\infty$ must be stationary.

\section{Short-time existence and parabolicity}
\setcounter{equation}{0}

To obtain short time existence, we consider the Anomaly flows as an evolution equation. Let $V$ be a smooth vector bundle over a compact manifold $X$, and consider the equation,
\bea
\label{evolution}
\p_t \psi={\cal E}(\psi)
\eea
where ${\cal E}(\psi)$ is a non-linear differential operator of order $2$, acting on the sections $\psi$ of $V$. If the eigenvalues of the symbol $\sigma(\delta{\cal E}(\psi))(x,\xi)$ of the linearization $\delta{\cal E}$ of ${\cal E}$
have strictly positive real parts for $\xi\not=0$, $(x,\xi)\in T_*(X)$, then the equation is parabolic, and the evolution equation with initial data $\psi$ admits a solution for short-time. 
More generally, we have the following version of the Nash-Moser theorem,
as formulated by Hamilton \cite{H}, and applied by him to show the existence of short-time solution for the Ricci flow:

\begin{lemma}
Let $L: C^\infty(V)\to C^\infty(W)$  be a linear differential operator of order $1$ with values in another vector bundle $W$. Assume that

\smallskip
{\rm (a)} The composition $Q(\Psi)=L(\Psi) {\cal E}(\Psi)$ is a differential operator of order at most $1$;

{\rm (b)} The symbol $\sigma(\delta {\cal E}(\Psi))(x,\xi)$ has eigenvalues with strictly positive real parts when restricted to the kernel of the symbol $\sigma(\delta L(\Psi))(x,\xi)$.

\smallskip
Then the initial value problem (\ref{evolution}) admits a unique solution for short time.

\end{lemma}

\subsection{Proof of Theorem \ref{Th2}}

We consider first the notationally simpler case of the Anomaly flow (\ref{anomalyX}) on $X$. In this case, the bundle $V$ is the bundle $\Lambda^{2,2}(X)$ of $(2,2)$-forms, the sections $\psi$ are the $(2,2)$-forms $\Psi$, and ${\cal E}(\psi)$
is given by the right hand side of (\ref{anomalyX}). The linearization of $i\ddb\o$ follows readily from the equation (\ref{variations1}),
\bea
\label{linearizationo}
\left(i\ddb \delta\omega\right)
=
i\ddb \left({1\over 2 \|\Omega\|_\o}
\left(\<\star\delta\Psi,\o\>\o
-\star\delta\Psi\right)\right)
=
i\ddb (\tilde\star\delta\Psi),
\eea
in the notation of (\ref{tildestar}). Next, we determine the linearization of the curvature terms in the Anomaly flow. The variation of the curvature $F$ of a unitary Chern connection 
under a variation $\delta H$ of the Hermitian metric is given by (see e.g. \cite{Siu}),
\bea
\delta F=  \bar\p \p^H(H^{-1}\delta H),
\eea
where $\p^H$ denotes the covariant exterior derivative in the unbarred directions. In particular,
\bea
\delta {\rm Tr}(F\wedge F)=2\,{\rm Tr}(F\wedge  \bar\p\p^H(H^{-1}\delta H)).
\eea
In view of the Bianchi identity, $d^HF=0$, this can be rewritten as
\bea
\delta {\rm Tr}(F\wedge F)
=
-2\, \ddb \,{\rm Tr}(F H^{-1}\delta H).
\eea

Specializing to the case where the vector bundle is $T^{1,0}(X)$, we obtain
\bea
\delta
{\rm Tr}(R\wedge R)
&=&
-2\ddb \left(  R_{\bar kj}{}^{\alpha\bar\beta}\delta g_{\bar\beta\alpha} dz^j\wedge d\bar z^k\right).
\nonumber
\eea
It follows that
\bea\label{linearizationRR}
\delta {\rm Tr} (R \wedge R) &=& 2 i \ddb \left(R_{\bar k j} {}^{\alpha\bar\beta} \, (\delta \o)_{\bar\beta\alpha} dz^j\wedge d\bar z^k\right) \nonumber\\
&=& 2i \ddb \left(Rm(\delta\o)\right)=
2i\ddb \left(Rm(\tilde\star\delta\Psi)\right)
\eea
where we view the Riemann curvature tensor as an operator on $(1,1)$-forms, as explained in \S 2.
Combining the formulas (\ref{linearizationo}) and (\ref{linearizationRR}) gives the linearization of ${\cal E}(\Psi)$ in the case of the Anomaly flow on $X$,
\bea
\delta{\cal E}(\delta\Psi)
=
\tilde\Delta (\delta\Psi).
\eea 

We apply Hamilton's version of the Nash-Moser implicit function theorem with the choice $L=d$ on $(2,2)$-forms. Because the right hand side ${\cal E}$ of the Anomaly flow is a closed $(2,2)$-form, we do have $L{\cal E}=0$. The condition of ellipticity of the operator $\tilde\Delta$ restricted to the space of closed $(2,2)$-forms, as formulated in \S 2, is precisely the condition which allows Hamilton's version of the Nash-Moser implicit function theorem to apply. Thus the existence of solutions to the equation (\ref{anomalyX}) for short-time follows.

\bigskip
The case of the Anomaly flow on $(X,E)$ can be treated in the same manner. We view the flow as of the form
\bea
\p_t\Psi={\cal E}(\Psi,H),
\qquad
\p_t H={\cal F}(\Psi,H)
\eea
with the pair $(\Psi,H)$ given by sections of the direct sum of the bundle of $(2,2)$-forms with the bundle of Hermitian quadratic forms on $E$. Clearly ${\cal E}$ and ${\cal F}$ are non-linear differential operators of order $2$.
Applying the formulas for variations of curvature to the bundles $T^{1,0}$ and $E$, we find
\bea
\delta {\cal E}
&=&
\ddb \bigg  (i \delta \o+ {\alpha' \over 2} \left({\rm Tr}(R g^{-1}\delta g)-{\rm Tr}(F H^{-1}\delta H)\right)\bigg)
\nonumber\\
\delta{\cal F}
&=&
-H\Lambda \bar\p\p^H(H^{-1}\delta H)
-
H \delta g^{j\bar k}F_{\bar kj}- \delta H\, \Lambda F,
\eea
where $F_{\bar kj}$ is viewed as a $(1,1)$-form with valued in the bundle of endomorphisms of $E$. It follows that the symbols of the linearization of ${\cal E}$ and ${\cal F}$ are given by
\bea
&&
\sigma(\delta{\cal E}):
\
(\delta\Psi,\delta H)
\to
\sigma(\delta{\cal E}_X)(\delta\Psi)
-
{\alpha' \over 2} \xi \wedge \bar\xi \wedge({\rm Tr}(FH^{-1}\delta H))
\nonumber\\
&&
\sigma(\delta{\cal F}):
\
(\delta\Psi,\delta H)
\to
|\xi|^2 \delta H
\eea
where we have temporarily denoted by ${\cal E}_X$ the expression on the right hand side of the Anomaly flow on $X$. 

We choose the operator $L$ of the Nash-Moser implicit function theorem as 
\begin{eqnarray*}
L(\Psi,H)=(d\Psi,0).
\end{eqnarray*} 
The right hand side ${\cal E}(\Psi,H)$ is again always a closed form, so we do have $L({\cal E},{\cal F})=0$. Furthermore,
the above formulas show that the symbol of the combined system $({\cal E},{\cal F})$ is a triangular block matrix, with the blocks on the diagonal given by $\sigma(\delta{\cal E}_X)$ acting on $\delta\Psi$ and $\sigma(\delta{\cal F})$ acting on $\delta H$. The first block has already been shown to have eigenvalues with positive real parts when restricted to the kernel of $d$, while the second is manifestly strictly positive. So the short-time existence of solutions to the Anomaly flow on $(X,E)$ follows, and the proof of Theorem \ref{Th2} is complete.

\subsection{Discussion of the parabolicity condition}

The linearization operator $\tilde\Delta$ and its ellipticity restricted to the space of $(2,2)$ closed forms is of considerable importance, since $\tilde\Delta$ controls the evolution of derivative quantities of $\Psi$ and $\o$ such as the curvature.

\smallskip
First we observe that the ellipticity condition can also be reformulated in terms of operators on $(1,1)$-forms. It suffices to set $\delta\Psi=\star\delta T$, for $(1,1)$-forms $T$. Then the ellipticity of $\tilde\Delta$ on the space of closed $(2,2)$-forms is equivalent to the ellipticity of the operator
\bea
\star\delta T\to \tilde \Delta(\star \delta T)
\eea
restricted to $(1,1)$-forms $\delta T$ satisfying $d^\dagger \delta T=0$, defined similarly in terms of the eigenvalues of its symbol, restricted to the kernel of the symbol of $d^\dagger$. This follows at once from the well-known formula $d^\dagger=-\star d\,\star$.

\smallskip
Next, it is instructive to work out the ellipticity condition of the operator $\tilde\Delta$ restricted to the space of $(2,2)$-forms more explicitly. The symbol $\sigma_{\tilde\Delta}$ of $\tilde\Delta$ is given by
\bea
\sigma_{\tilde\Delta}(\xi):
\ \delta\Psi\ \rightarrow
\
i\xi\wedge\bar\xi\wedge (\tilde\star\delta\Psi
-{\alpha' \over 2}Rm(\tilde\star\delta\Psi))
\eea
where $\xi\in T^*(X)$. The following lemma is useful:

\begin{lemma}
\label{lemma}
For all $\delta\Psi$ in the kernel of the symbol of the exterior derivative $d$ on $(2,2)$-forms, we have
\bea
\label{Id}
i\xi\wedge\bar \xi \wedge \tilde\star \delta\Psi
=
{1\over 2\|\Omega\|_\o}|\xi|^2\,\delta\Psi.
\eea
\end{lemma}

\medskip
{\it Proof.} Without loss of generality, we can choose coordinates so that $g_{\bar jk}=\delta_{jk}$, and by an additional orthogonal rotation if necessary, that $\xi=(\xi_1, 0, 0)$. Let $\Phi$ denote the left-hand side of (\ref{Id}). Using the coordinate expression (\ref{variations}) and the convention (\ref{n-1_convention}) for components of a $(2,2)$-form, we find
%In the component notation for $(2,2)$-forms in dimension $3$, its components are given by $\Phi^{\ell\bar m}=\e^{\ell rs}\e^{\bar m\bar p\bar q}\Phi_{\bar p r\bar qs}$. We find
\bea
\label{Phi}
&&
\Phi^{k\bar 1}=0, \qquad k=1,2,3
\nonumber\\
&&
\Phi^{k\bar k}
=
{1\over 2 \|\Omega\|_\o}
|\xi|^2(\delta\Psi^{1\bar 1}+\delta\Psi^{k\bar k}),
\qquad k=2,3
\nonumber
\\
&&
\Phi^{2\bar 3}={1\over 2 \|\Omega\|_\o}|\xi|^2 \delta\Psi^{2\bar 3}.
\eea
Now the kernel of the exterior derivative $d$ on $(2,2)$-forms is given by forms $\delta\Psi$ satisfying
\bea
\xi_j\delta\Psi^{j\bar k}=0\qquad k=1,2,3.
\eea
In the given coordinate system, this reduces to $\delta\Psi^{1\bar k}=0$ for $k=1,2,3$. Using these relations, we can rewrite the above identities as $\Phi= {1 \over 2} \|\Omega\|_\o^{-1}|\xi|^2\delta\Psi$. The lemma is proved.

\bigskip

This allows us to identify immediately an important and quite general situation where the ellipticity condition, and hence the existence of short-time solutions, holds:

\begin{proposition} \label{short-time-condition}
Consider the operator 
\bea
\delta\Psi
\to - i\xi\wedge\bar\xi \wedge {\alpha' \over 2} Rm(\<\star\delta\Psi,\o\>\o-\star\delta\Psi)
\eea
restricted to the kernel of the symbol of the exterior derivative $d$ on $(2,2)$-forms. If it has operator norm $<|\xi|^2$ for any $\xi\not=0$,
then the Anomaly flows with $\o$ as initial metric admit smooth solutions for at least a short time.
\end{proposition}

\section{Remarks}
\setcounter{equation}{0}

We conclude with a few remarks.

\medskip

(1) {\bf The balanced condition}

\medskip
One of the major challenges in Strominger systems is that, even with the metric $H$ on the vector bundle $E$ fixed, it is a system in the metric $\o$, in the sense that both the anomaly equation and the balanced-like condition have to be satisfied. 
One natural approach is to try and solve the anomaly equation with a particular ansatz which guarantees that the metric is balanced. One possible ansatz is the very general one proposed by Fu-Wang-Wu \cite{FWW}, Tosatti-Weinkove \cite{TW1} and Popovici \cite{P}
\bea
\label{TW}
\o^2=\o_0^2+i\ddb(u\tilde\o)
\eea
Here $\o_0$ is a balanced metric, and $\tilde\o$ is an arbitrary $(1,1)$-form, and the form $\o^2$ is required to be positive. Tosatti and Weinkove have also shown in \cite{TW1} how to find a single balanced metric $\o_0$.

Another ansatz is the one by Fu and Yau \cite{FY}, in the special case of Goldstein-Prokushkin manifolds discussed below (see eq. (\ref{FY}) below). But unlike in the K\"ahler case, where metrics can be represented by a potential, there does not appear to be a uniquely compelling ansatz for balanced metrics at the present time. Thus, a 
very attractive feature of the Anomaly flows is that they guarantee that the metrics be balanced, without appealing to any particular ansatz.

\

\medskip
(2) {\bf Toric fibrations over $K3$ surfaces}

\medskip
\par
The first non-perturbative solution of a Strominger system was found by Fu-Yau \cite{FY}, as a toric fibration over a $K3$ surface. More specifically, Goldstein and Prokushkin \cite{GP} had shown how to construct a toric fibration $\pi:X\to S$, given a Calabi-Yau surface $(S,\o_S)$ with Ricci-flat K\"ahler metric $\o_S = i \, (g_S){}_{\bar{k} j} dz^j \wedge d \bar{z}^k$, and two anti-self-dual $(1,1)$-forms $\kappa_1,\kappa_2\in 2\pi H^2(S,{\bf Z})$. Furthermore, there is a $(1,0)$-form $\theta$ on $X$ so that $\p \theta =0$, $\bar{\p} \theta = \pi^* (\kappa_1 + i \kappa_2)$, and if $\Omega_S$ is a non-vanishing holomorphic $(2,0)$-form on $S$, then the non-vanishing holomorphic $(3,0)$-form on $X$ is given by $\Omega=\theta\wedge\pi^*(\Omega_S)$. The $(1,1)$-form $\o_0=\pi^*(\o_S)+i\theta\wedge\bar\theta$ is a balanced metric, $d\o_0^2=0$ and also satisfies $||\Omega||_{\omega_0} =1$. This implies that $d\left(||\Omega||_{\omega_0} \o_0^2 \right)=0$.

Under suitable cohomological conditions on the class $\kappa_1$ and $\kappa_2$, Fu and Yau found a solution of the Strominger system under the following ansatz,
\bea
\label{FY}
\o_u=\pi^*(e^u\o_S)+i\theta\wedge \bar\theta
\eea
where $u$ is a scalar function on $S$. Metrics of the form $\o_u$ are automatically balanced and satisfy $\| \Omega \|_{\o_u} = e^{-u}$. Here we observe that another hint that the Anomaly flow is a natural flow, is that it preserves the Fu-Yau ansatz (\ref{FY}). Indeed, under the cohomological conditions mentioned previously, Fu and Yau have shown that
\bea\label{FuYauequ}
i\ddb\o_u
-
{\alpha' \over 4}({\rm Tr}(R_u\wedge R_u)-{\rm Tr}(F\wedge F))
=
i\ddb(e^u \o_S - \alpha' e^{-u} \rho)
-
{\alpha' \over 2} \, \ddb u\wedge\, \ddb u+
\tilde{\mu} {\o_S^2\over 2!}
\nonumber
\eea
for some smooth function $\tilde{\mu}:S \rightarrow {\bf R}$ and smooth real $(1,1)$ form $\rho$ defined on $S$ and given by
\be \label{defn_rho}
\rho = - {i \over 2} \, {\rm Tr} (\bar\p B\wedge \p B^* \, g_S^{-1}),
\ee
where $B$ is defined on $S$ and only depends on $(S,\o_S)$ and $\kappa_1,\kappa_2$. Thus $\p_t(\|\Omega\|_{\o_u}\o_u^2)$ is the pull-back of a $(2,2)$-form on $S$. Since
\bea
\|\Omega\|_{\o_u}\o_u^2
=
\o_0^2+(e^u-1)\o_S^2
\eea
we see that an evolution of $\|\Omega\|_{\o_u}\o_u^2$ by a term proportional to $\o_S^2$ is just an evolution of the conformal factor $u$, and $\o_u$ still satisfies the same ansatz. In fact, the Anomaly flow is immediately seen to be equivalent to the equation
\bea
\label{FYanomaly}
\p_tu
=
{e^{-u} \over 2} \bigg(\Delta e^u + \alpha' \sigma_2( i \ddb u) -2 \alpha' {i \ddb (e^{-u} \rho) \over \o_S^2} +\tilde{\mu} \bigg),
\eea
where the Laplacian and $\sigma_2(i \ddb u) \o_S^2 = i \ddb u \wedge i \ddb u$ are with respect to the metric $\o_S$.
This equation may be of interest in its own right. Its parabolicity is
equivalent to the ellipticity of the right hand side, which is the ellipticity condition imposed by Fu and Yau (\cite{FY}, eqs. (7.3) and (8.1)). Explicitly, 
denote
\begin{eqnarray*}
E(u) &=&i\ddb\o_u
-
{\alpha' \over 4} \left({\rm Tr}(R_u\wedge R_u)-{\rm Tr}(F\wedge F)\right)\\
&=& i\ddb(e^u \o_S - \alpha' e^{-u} \rho)
+
{\alpha' \over 2} \,i\ddb u\wedge\,i\ddb u+
\tilde{\mu} \, {\o_S^2 \over 2!}.
\end{eqnarray*}
Its symbol is given by
\bea
 \sigma(\delta E)\ : \ \delta u \rightarrow i\xi\wedge \bar{\xi} \wedge (\delta u )\left( e^u  \omega_S + \alpha' e^{-u} \rho + \alpha' i\partial\bar\partial u \right).
\eea
The ellipticity condition reduces to the following condition on $u$
\bea \label{FY_ellip}
e^u  \omega_S + \alpha' e^{-u} \rho + \alpha' i\partial\bar\partial u  >0.
\eea
We now compare this condition with Proposition \ref{short-time-condition}. The condition
\be \label{curv_small}
|\alpha' Rm| \ll 1,
\ee
implies short-time existence for the Anomaly flow by Proposition \ref{short-time-condition}. In \cite{FY}, Fu and Yau computed the curvature of the metric $\o_u$ (\ref{FY}). Fixing a point $p \in X$, they constructed a frame of holomorphic vector fields such that at $p$,
\be
g_u = \left( \begin{array}{cc}
e^u g_S & 0\\
0 & 1\end{array}\right) \ \ \ Rm=\left( \begin{array}{cc}
R_{11} & R_{12}\\
R_{21} & R_{22}\end{array}\right)
\ee
where the entries $R_{jk}$ are given by
\bea
R_{11}&=& R_S - \ddb u \, I +e^{-u} \bar\p B\wedge \p B^* \,g_S ^{-1}\\
R_{12}&=& -\na\bar\p B+\p u\wedge\bar\p B\\
R_{21}&=&\bar\p( e^{-u} \p B^* \,g_S^{-1})\\
R_{22}&=& e^{-u} (\p B^*\,g_S^{-1})\wedge \bar\p B.
\eea
Here $B=(\phi_1,\phi_2)^T$ is a column vector of locally defined functions $\phi_i$ on $S$, and $\bar{\p} B$ is globally defined on $S$. Using the definition (\ref{defn_rho}) of $\rho$ and the fact that ${\rm Tr} R_S=0$, we see that $|\alpha' Rm|_{\o_u} \ll 1$ implies
\be \label{FY_key_est}
|2 \alpha' e^{-2u} \rho|_{\o_S} \ll 1, \ \ |2 \alpha' e^{-u} \ddb u|_{\o_S} \ll 1.
\ee
It is clear that (\ref{FY_key_est}) implies (\ref{FY_ellip}).
\medskip

\bigskip

{\bf Update:}

\medskip

We note that there has been very recently significant progress on the Anomaly flows in several directions. In particular, an explicit formula for $\p_tg_{\bar kj}$ has now been obtained \cite{PPZ1} which overcomes the difficulty of the flow being originally formulated as a flow of $(2,2)$-forms. Furthermore, the convergence of the Anomaly flow on Goldstein-Prokushkin fibrations has been established in \cite{PPZ2}, for both $\alpha'>0$ and $\alpha'<0$, thus unifying the results of Fu and Yau \cite{FY, FY2}. These constitute strong evidence that the Anomaly flows should provide an efficient approach to Strominger systems.

\

%\newpage

\bigskip
Department of Mathematics, Columbia University, New York, NY 10027, USA

\smallskip

phong@math.columbia.edu

\bigskip
Department of Mathematics, Columbia University, New York, NY 10027, USA

\smallskip
picard@math.columbia.edu

\bigskip
Department of Mathematics, University of California, Irvine, CA 92697, USA

\smallskip
xiangwen@math.uci.edu


\begin{thebibliography}{99}

{\small

\bibitem{AG1} Andreas, B., and Garcia-Fernandez, M., {\it Heterotic non-Kahler geometries via
polystable bundles on Calabi-Yau threefolds}, Journal of Geometry and Physics, Vol 62 (2012), no. 2, 183-188.

\bibitem{AG2} Andreas, B., and Garcia-Fernandez, M., {\it Solutions of the Strominger system via stable bundles on Calabi-Yau threefolds}, Communications in Mathematical Physics, 315 (2012), 153-168.

\bibitem{BBFT} Becker, K., Becker, M., Fu, J.X., Tseng, L.S., and Yau, S.T., {\em Anomaly cancellation and smooth non-Kahler solutions in heterotic string theory}, Nuclear Physics B, 751 (2006) no. 1, 108-128.

\bibitem{CI} Carlevaro, L., and Israel, D., {\em Heterotic resolved conifolds with torsion, from
supergravity to CFT}, Journal of High Energy Physics, (2010) no.1, 1-57.

\bibitem{D} Donaldson, S., {\em Infinite determinants, stable bundles, and curvature}, Duke Math. J. 54 (1987) 231-247.

\bibitem{DRS} Dasgupta, K., Rajesh, G., and Sethi, S.S., {\em M-theory, orientifolds and G-flux}, Journal of High Energy Physics, no. 08 (1999).

\bibitem{Fe1} Fei, T., {\em A construction of non-K\"ahler Calabi-Yau manifolds and new solutions to the Strominger system}, Adv. Math. 302, 529-550 (2016).

\bibitem{Fe2} Fei, T., {\em Some Torsional Local Models of Heterotic Strings}, preprint, arXiv:1508.05566.

\bibitem{FeY} Fei, T. and Yau, S.T., {\em Invariant solutions to the Strominger system on complex Lie groups and their quotients}, Comm. Math. Phys., Vol 338, Number 3 (2015), 1-13.

\bibitem{FIUV1} Fernandez, M., Ivanov, S., Ugarte, L., and Vassilev, D., {\em Non-Kahler heterotic string solutions with non-zero fluxes and non-constant dilaton}, Journal of High Energy Physics, Vol 6, (2014) 1-23.

\bibitem{FIUV2} Fernandez, M., Ivanov, S., Ugarte, L., and Villacampa, R., {\em Non-Kahler heterotic string compactifications with non-zero fluxes and constant dilaton}, Comm. Math. Phys. 288 (2009), 677-697.

\bibitem{FTY} Fu, J.X., Tseng, L.S., and Yau, S.T., {\em Local heterotic torsional models}, Communications in Mathematical Physics, 289 (2009), 1151-1169.

\bibitem{FWW} Fu, J.X., Wang, Z.Z. and Wu, D.M., {\em Form-type equations on K\"ahler manifolds of nonnegative orthogonal bisectional curvature}, Calc. Var. Partial Differential Equations 52 (2015), no. 1-2, 327-344.

\bibitem{FY} Fu, J.X. and Yau, S.T. {\it The theory of superstring with flux on non-Kahler manifolds and the complex Monge-Ampere equation}, J. Differential Geom, Vol 78, Number 3 (2008), 369-428.

\bibitem{FY2} Fu, J.X. and Yau, S.T., {\em A Monge-Amp\`ere type equation motivated by string theory}, Comm. Anal. Geom., Vol 15, Number 1 (2007), 29-76.


\bibitem{GP} Goldstein, E. and Prokushkin, S., {\it Geometric model for complex non-K\"ahler manifolds with SU(3) structure}, Comm. Math. Phys. 251(2004), no. 1, 65-78.

\bibitem{Gr} Grantcharov, G., {\em Geometry of compact complex homogeneous spaces with vanishing first Chern class}, Advances in Mathematics, 226 (2011), 3136-3159.

\bibitem{H} Hamilton, R.S, {\it Three-Manifolds with Positive Ricci Curvature}, J. Differential Geom. 17 (1982), 255-306.

\bibitem{Hull1} Hull, C., {\em Superstring Compactifications with Torsion and Space-Time Supersymmetry}, Proceedings of the First Torino Meeting on Superunification and Extra Dimensions, edited by R.D'Auria and P.Fre, World Scientific, Singapore, (1986).

\bibitem{Hull2} Hull, C., {\em Compactifications of the heterotic superstring}, Physics Letters B 178 (1986), no. 4, 357-365.

\bibitem{LY} Li, J. and Yau, S.T., {\em The existence of supersymmetric string theorey with torsion}, J. Differential Geom. 70 (2005), no.1, 143-181.

\bibitem{M} Michelsohn, M.L., {\it On the existence of special metrics in complex geometry}, Acta Math. 149 (1982), no. 3-4, 261-295.

\bibitem{OUV} Otal, A., Ugarte, L., Villacampa, R. {\em Invariant solutions to the Strominger system and the heterotic equations of motion on solvmanifolds}, preprint arXiv:1604.02851.

\bibitem{PPZ} Phong, D.H., Picard, S. and Zhang, X.W., {\em On estimates for the Fu-Yau generalization of a Strominger system}, to appear in J. Reine Angew. Math. (arXiv:1507.08193)

\bibitem{PPZ0} Phong, D.H., Picard, S. and Zhang, X.W., {\em The Fu-Yau equation with negative slope parameter}, to appear in Invent. Math. (arXiv:1602.08838).

\bibitem{PPZ1} Phong, D.H., Picard, S. and Zhang, X.W., {\em Anomaly flows}, arXiv:1610.02739.

\bibitem{PPZ2} Phong, D.H., Picard, S. and Zhang, X.W., {\em The Anomaly flow and the Fu-Yau equation}, arXiv:1610.02740.

\bibitem{P} Popovici, D., {\em Aeppli cohomology classes associated with Gauduchon metrics on compact complex manifolds}. Bulletin de la Soci\'et\'e Math\'ematique de France, 143(4), 763-800.

\bibitem{Siu} Siu, Y.-T., {\it Lectures on Hermitian-Einstein metrics for stable bundles and K\"ahler-
Einstein metrics}, DMV Seminar, 8. Birkh\"auser Verlag, Basel, 1987.

\bibitem{S} Strominger, A., {\it Superstrings with Torsion}, Nuclear Physics B 274 (1986), no.2, 253-284.

\bibitem{TW1} Tosatti, V. and B. Weinkove, {\it The Monge-Amp\`ere equation for $(n-1)$-plurisubharmonic functions on a compact K\"ahler manifold}, 	J. Amer. Math. Soc. 30 (2017), no.2, 311-346..


\bibitem{UV} Ugarte, L., Villacampa, R., {\em Non-nilpotent complex geometry of nilmanifolds and heterotic supersymmetry}, Asian Jour. Math. 18 (2014), no 2, 229-246.

\bibitem{UY} Uhlenbeck, K. and S.T. Yau,
{\em On the existence of Hermitian Yang-Mills connections on stable vector bundles}, Comm. Pure Appl. Math. 39 (1986) 257-293.

\bibitem{Y} Yau, S.T., {\it On the Ricci curvature of a compact K\"ahler manifold and the complex Monge-Amp\`ere equation}, I, Comm. Pure Appl. Math. 31 (1978), no.3, 339-411.





}

\end{thebibliography}
\end{document}